# Asymptotics for sums of functions of primes

VICTOR VOLFSON

ABSTRACT This work gives a general approach to the determination of the asymptotic behavior of the sums of functions of primes based on the distribution of primes. It refines the estimate of the remainder term of the asymptotic expansion of the sums of functions of primes. Also, the necessary and sufficient conditions for the existence of these asymptotics are proved in the paper.

1. INTRODUCTION

Using the method of complex integration the following asymptotic of the number of primes not exceeding the natural value $n$ is shown in [1]:

$$\pi(n) = \sum_{p \leq n} 1 = \int_2^n \frac{du}{\log(u)} + O\left(\frac{n}{e^{c\sqrt{\log(n)}}}\right), \qquad (1.1)$$

where $c > 0$ is the constant.

Based on (1.1), we can obtain a more crude asymptotic estimate of the number of primes not exceeding the natural value $n$. However, it is sometimes convenient to use.

Integrating by parts, we get:

$$\int_2^n \frac{du}{\log(u)} = \frac{u}{\log(u)}\Big|_2^n - \int_2^n u\, d\left(\frac{1}{\log(u)}\right) = \frac{n}{\log(n)} - \frac{2}{\log 2} + \int_2^n \frac{du}{\log^2 u} = \frac{n}{\log(n)} + O\left(\frac{n}{\log^2 n}\right)$$

Substitute this expression in (1.1) and obtain the asymptotic:

$$\pi(n) = \frac{n}{\log(n)} + O\left(\frac{n}{\log^2 n}\right) + O\left(\frac{n}{e^{c\sqrt{\log(n)}}}\right) = \frac{n}{\log(n)} + O\left(\frac{n}{\log^2 n}\right). \qquad (1.2)$$

---

Expression (1.1) for the asymptotic distribution of prime numbers was first obtained in 1896 by Hadamard and Valais-Poussin [2] based on an estimate of the boundary of the region of absence of nontrivial zeros of the Riemann zeta function:

$$\sigma \geq 1 - \frac{c}{\log T}, \qquad (1.3)$$

where $T \geq 10$.

Vinogradov [3] (in 1958) proved a new bound for the boundary of the region of absence of nontrivial zeros of the Riemann zeta function:

$$\sigma \geq 1 - \frac{c}{\log T^{2/3-\xi}}. \qquad (1.4)$$

Based on (1.4), a new asymptotic estimate for the remainder in formula (1.1) was obtained:

$$O\left(\frac{n}{e^{c(\log(n))^{0,6-\xi}}}\right). \qquad (1.5)$$

The result (1.5) is practically not improved to this day. Meanwhile, on the basis of the Riemann hypothesis [4], which has not yet been proved, all non-trivial zeros of the zeta function are on the straight line $\sigma = 1/2$.

The following asymptotic estimate for the number of primes not exceeding a natural value $n$ is carried out if the Riemann conjecture is true:

$$\pi(n) = \int_2^n \frac{du}{\log(u)} + O(n^{1/2} \log(n)). \qquad (1.6)$$

Based on the indicated asymptotics for the number of primes not exceeding a natural value $n$, a general formula will be obtained for the asymptotic estimate of the sums and products of functions of primes in the work.

It is logical that based on the asymptotic law of distribution of prime numbers, asymptotics of the sums of functions of prime numbers are obtained.

Previously, there was no such general approach; therefore, asymptotic estimates for the sums of functions of prime numbers were very laborious, not systematized, and were obtained from other considerations [5].

Now, in the theory of primes, the possibility of a systematic approach has appeared - first, to state the asymptotic law of distribution of primes and the Riemann conjecture, and then, based on them, to derive a general formula for the asymptotic estimation of sums of primes, and then, as an example of using the general formula, to derive individual formulas as done in this work.

The following general formula was obtained in [6]. However, this formula contains only the main term of the expansion. Therefore, the purpose of this work is to obtain the remainder of the expansion and consider various examples of determining the asymptotics of the sums of functions of prime numbers using the resulting general formula, which will be done in Chapter 2 of this work.

Necessary and sufficient conditions for the existence of asymptotics will be proved in Chapter 3.

## 2. ASYMTOTIC OF SUMS OF FUNCTIONS OF PRIME NUMBERS

First, we use formula (1.2) to determine the asymptotics of the sums of functions of prime numbers.

Assertion 2.1

If $f'(t)$ exists and is continuous, then:

$$\sum_{p \leq n} f(p) = \frac{nf(n)}{\log(n)} + O(\frac{n|f(n)|}{\log^2(n)}) - \int_2^n \frac{tf'(t)dt}{\log(t)} + O(\int_2^n \frac{t|f'(t)|dt}{\log^2(t)})$$

Proof

Let $a_k = 1$, if $k$ is a prime number and $a_k = 0$ otherwise. Let us denote $A(n) = \sum_{k=1}^{n} a_k = \pi(n)$.

If $f'(t)$ exists and is continuous, then we obtain, using the Abel's formula for the sum $\sum_{k=1}^{n} a_k f(k)$:

$$\sum_{p \leq n} f(p) = \sum_{k=1}^{n} a_k f(k) = A(n) f(n) - \int_{1}^{n} A(t) f'(t) dt, \qquad (2.1)$$

where $A(n) = \sum_{k=1}^{n} a_k = \pi(n)$.

In accordance with the formula (1.2) the value:

$$A(n) = \pi(n) = \frac{n}{\log(n)} + O(n / \log^2(n)). \qquad (2.2)$$

Substituting (2.2) into (2.1) and get:

$$\sum_{p \leq n} f(p) = \frac{nf(n)}{\log(n)} + O(\frac{n |f(n)|}{\log^2(n)}) - \int_{2}^{n} \frac{tf'(t)dt}{\log(t)} + O(\int_{2}^{n} \frac{t|f'(t)|dt}{\log^2(t)}), \qquad (2.3)$$

which corresponds to Assertion 2.1.

Let us give examples of using formula (2.3):

1. $\sum_{p \leq n} 1 = \frac{n}{\log(n)} + O(\frac{n}{\log^2(n)})$, which corresponds to formula (1.2).

2. $\sum_{p \leq n} \log(p) = n + O(n / \log(n))$. $\qquad (2.4)$

3. $\sum_{p \leq n} 1/p = \log\log(n) + O(1)$. $\qquad (2.5)$

Now we use a more precise asymptotic estimate (1.1).

Assertion 2.2

Suppose that the function $f$ is monotone and has a continuous derivative on the interval $[2, n]$, then:

$$\sum_{p \leq n} f(p) = \int_{2}^{n} \frac{f(t)dt}{\log(t)} + O(\frac{|f(n)|n}{e^{c\log^{1/2}(n)}}) + O(\int_{2}^{n} \frac{t|f'(t)|dt}{e^{c\log^{1/2}(n)}}).$$

Proof

Let $a_k = 1$, if $k$ is a prime number and otherwise- $a_k = 0$.

We denote: $A(n) = \pi(n) = \sum_{k=1}^{n} a_k$.

Having in mind (1.1):

$$A(n) = \pi(n) = \int_{2}^{n} \frac{dt}{\log(t)} + O(\frac{n}{e^{c\log^{1/2}(n)}}).  \qquad (2.6)$$

If the function $f$ is monotone and it has a continuous derivative on the interval $[2, n]$, then, based on the Abel's summation formula, we can write:

$$\sum_{p \leq n} f(p) = \sum_{k=1}^{n} a_k f(k) = A(n)f(n) - \int_{2}^{n} A(t) f'(t) dt. \qquad (2.7)$$

Substituting (2.6) into (2.7) and get:

$$\sum_{p \leq n} f(p) = f(n) \int_{2}^{n} \frac{du}{\log(u)} + O(\frac{|f(n)|n}{e^{c\log^{1/2}(n)}}) - \int_{2}^{n} (\int_{2}^{t} \frac{du}{\log(u)}) f'(t) dt + O(\int_{2}^{n} \frac{t|f'(t)|dt}{e^{c\log^{1/2}(n)}}). \qquad (2.8)$$

We use the method of integration by parts:

$$\int_{2}^{n} (\int_{2}^{t} \frac{du}{\log(u)}) f'(t) dt = f(n) \int_{2}^{n} \frac{du}{\log(u)} - \int_{2}^{n} \frac{f(t) dt}{\log(t)}. \qquad (2.9)$$

Substituting (2.9) into (2.8) and get:

$$\sum_{p \leq n} f(p) = \int_{2}^{n} \frac{f(t) dt}{\log(t)} + O(\frac{|f(n)|n}{e^{c\log^{1/2}(n)}}) + O(\int_{2}^{n} \frac{t|f'(t)|dt}{e^{c\log^{1/2}(n)}}), \qquad (2.10)$$

which corresponds to Assertion 2.2.

Let's take a look at examples of using (2.10):

1. $\sum_{p \leq n} 1 = \int_{2}^{n} \frac{dt}{\log(t)} + O(\frac{n}{e^{c\log^{1/2}(n)}})$, which corresponds to (1.1).

2. $\sum_{p \leq n} \log(p) = \int_{2}^{n} \frac{\log(t) dt}{\log(t)} + O(\frac{\log(n)n}{e^{c\log^{1/2}(n)}}) + O(\int_{2}^{n} \frac{dt}{e^{c\log^{1/2}(n)}}) = n + O(\frac{\log(n)n}{e^{c\log^{1/2}(n)}}). \qquad (2.11)$

Let us compare (2.11) and (2.4).

1. $\sum_{p \le n} \frac{\log(p)}{p} = \int_2^n \frac{\log(t)dt}{t\log(t)} + O(\frac{\log(n)n}{ne^{c\log^{1/2}(n)}}) + O(\int_2^n \frac{tdt}{t^2 e^{c\log^{1/2}(n)}}) + O(\int_2^n \frac{t\log(t)dt}{t^2 e^{c\log^{1/2}(n)}}) = \log(n) + O(1)$

2. $\sum_{p \le n} p^\alpha = \int_2^n \frac{t^\alpha dt}{\log(t)} + O(\frac{n^{\alpha+1}}{e^{c\log 1/2(n)}}) + O(\int_2^n \frac{t^\alpha dt}{e^{c\log 1/2(t)}}) = \int_2^n \frac{t^\alpha dt}{\log(t)} + O(\frac{n^{\alpha+1}}{e^{c\log 1/2(n)}})$,

where $\alpha > -1$.

Assertion 2.3

Let the function $f$ is monotone and it has a continuous derivative on the interval $[2,n]$, then if the Riemann conjecture is true, the best estimate of the remainder of the asymptotic is fulfilled:

$$\sum_{p \le n} f(p) = \int_2^n \frac{f(t)dt}{\log(t)} + O(|f(n)| n^{1/2} \log(n)) + O(\int_2^n |f'(t)| t^{1/2} \log(t) dt)$$

Proof

The proof is similar to Assertion 2.2, only having in mind (1.6) in the Abel's formula the value $A(n)$ is taken as:

$$A(n) = \pi(n) = \int_2^n \frac{du}{\log(u)} + O(n^{1/2} \log(n)).$$

The remainder estimate in assertion 2.3 is the best. This can be seen by comparing the remainder terms in this assertion with s assertions 2.1 and 2.2.

Now consider an example of an asymptotic upper bound for the product of primes.

Let us denote $A = \prod_{p \le n} p$. Then:

$$\log A = \sum_{p \le n} \log(p). \tag{2.12}$$

Earlier (2.4) we defined the asymptotic of such a sum:

$$\sum_{p \le n} \log(p) = n + O(n/\log(n)). \tag{2.13}$$

If the Riemann conjecture is true, then, based on assertion 2.3, the asymptotic of this sum is:

$$\sum_{p \leq n} \log(p) = n + O(n^{1/2+\epsilon}). \qquad (2.14)$$

Based on (2.12) and (2.13), we obtain the following asymptotic upper bound for the product of primes:

$$\prod_{p \leq n} p \leq e^{n + c_1 n / \log(n)}, \qquad (2.15)$$

where constant $c_1 > 0$.

If the Riemann conjecture is true, then, having in mind (2.12) and (2.14), it holds the following asymptotic upper bound for the product of primes:

$$\prod_{p \leq n} p \leq e^{n + c_2 n^{1/2+\epsilon}}, \qquad (2.16)$$

where constant $c_2 > 0$.

Let us compare the asymptotic upper bounds (2.15) and (2.16). Let us compare assertions 2.1 and 2.2, i.e. formulas (2.3) and (2.10).

The principal terms of the asymptotics (2.3) and (2.10) must coincide, and the differences can only be in the remainder terms. Let's show it.

We transform the term in the formula (2.10):

$$\int_2^n \frac{f(t)dt}{\log(t)} = \frac{f(n)n}{\log(n)} - \frac{2f(2)}{\log(2)} - \int_2^n \frac{tf'(t)dt}{\log(t)} + \int_2^n \frac{f(t)dt}{\log(t)}. \qquad (2.17)$$

Substitute (2.17) into (2.10):

$$\sum_{p \leq n} f(p) = \frac{f(n)n}{\log(n)} + O(1) - \int_2^n \frac{tf'(t)dt}{\log(t)} + \int_2^n \frac{f(t)dt}{\log(t)} + O(\frac{|f(n)n|}{e^{c \log^{1/2}(n)}}) + O(\int_2^n \frac{t|f'(t)|dt}{e^{c \log^{1/2}(n)}}). (2.18)$$

Let us compare expression (2.18) with the expression for asymptotic (2.3). The first three terms are the same. The discrepancy is only in the residual terms. The remainder terms of expression (2.18) give a more accurate estimate.

# 3. NECESSERY AND SUFFICIENT CONDITIONS FOR THE EXISTENCE OF THE ASYMTOTICS

Let $a_k = 1$, if $k$ is a prime and $a_k = 0$, if $k$ is composite and let $b_k = \dfrac{1}{\ln k}, b_1 = 0$.

We denote: $A(n) = \sum_{k \leq n} a_k$ and $B(n) = \sum_{k \leq n} b_k$.

It is required that:

$$\lim_{n \to \infty} \frac{A(n)}{B(n)} = 1. \tag{3.1}$$

It follows from (3.1), that it is necessary to satisfy the indicated asymptotic:

$$\lim_{n \to \infty} \frac{\sum_{k=1}^{n} a_k f(k)}{\sum_{k=1}^{n} b_k f(k)} = 1. \tag{3.2}$$

Assertion 3.1

If the specified conditions are met:

1. The limit $\lim_{n \to \infty} \dfrac{\int_{1}^{n} B(t) f'(t) dt}{B(n) f(n)}$ is not equal to 1.

2. $f(x)$ is monotonous and $f'(x) \neq 0$.

3. The limit $\lim_{n \to \infty} \int_{1}^{n} B(t) f'(t) dt = \pm \infty$.

Then:

$$\lim_{n \to \infty} \frac{\sum_{k=1}^{n} a_k f(k)}{\sum_{k=1}^{n} b_k f(k)} = 1$$

Proof

Based on Abel's summation formula:

$$\sum_{k=1}^{n} a_k f(k) = A(n)f(n) - \int_1^n A(t)f'(t)dt, \quad \sum_{k=1}^{n} b_k f(k) = B(n)f(n) - \int_1^n B(t)f'(t)dt. \quad (3.3)$$

Having in mind (3.2) and (3.3) we get:

$$\frac{\sum_{k=1}^{n} a_k f(k)}{\sum_{k=1}^{n} b_k f(k)} = \frac{A(n)f(n) - \int_1^n A(t)f'(t)dt}{B(n)f(n) - \int_1^n B(t)f'(t)dt} = \frac{A(n)}{B(n)} \cdot \frac{1 - \dfrac{\int_1^n A(t)f'(t)dt}{A(n)f(n)}}{1 - \dfrac{\int_1^n B(t)f'(t)dt}{B(n)f(n)}}. \quad (3.4)$$

Suppose that the following conditions are met:

1. The limit $\lim\limits_{n \to \infty} \dfrac{\int_1^n B(t)f'(t)dt}{B(n)f(n)}$ is not equal to 1.

2. Let $f(x)$ is monotonous and $f'(x) \neq 0$.

3. The limit $\lim\limits_{n \to \infty} \int_1^n B(t)f'(t)dt = \pm\infty$.

Let us show that under these conditions:

$$\lim_{n \to \infty} \frac{\int_1^n A(t)f'(t)dt}{A(n)f(n)} = \lim_{n \to \infty} \frac{\int_1^n B(t)f'(t)dt}{B(n)f(n)}. \quad (3.5)$$

Indeed, in this case, using the L'Hôpital rule, we get:

$$\lim_{n \to \infty} \frac{\dfrac{\int_1^n A(t)f'(t)dt}{A(n)f(n)}}{\dfrac{\int_1^n B(t)f'(t)dt}{B(n)f(n)}} = \lim_{n \to \infty} \frac{B(n)}{A(n)} \cdot \lim_{n \to \infty} \frac{\int_1^n A(t)f'(t)dt}{\int_1^n B(t)f'(t)dt} = \lim_{n \to \infty} \frac{A(n)f'(n)}{B(n)f'(n)} = 1, \quad (3.6)$$

which corresponds to (3.5).

Then, based on (3.1) and (3.6), we obtain:

$$\lim_{n\to\infty} \frac{\sum_{k=1}^{n} a_k f(k)}{\sum_{k=1}^{n} b_k f(k)} = \lim_{n\to\infty} \frac{A(n)}{B(n)} \cdot \frac{1 - \frac{\int_1^n A(t)f'(t)dt}{A(n)f(n)}}{1 - \frac{\int_1^n B(t)f'(t)dt}{B(n)f(n)}} = 1,$$

which is consistent with the assertion.

Corollary 3.2

Conditions (1), (3) in assertion 3.1 correspond to:

1. The limit $\lim\limits_{n\to\infty} \dfrac{\int_2^n \frac{tf'(t)}{\log(t)}dt}{\frac{nf(n)}{\log(n)}}$ is not equal to 1.

3. The limit $\lim\limits_{n\to\infty} \int_2^n \dfrac{tf'(t)}{\log(t)}dt = \pm\infty$.

Proof

Based on the asymptotic law of prime numbers:

$$A(n) = B(n) = \frac{n}{\log(n)}(1 + o(1)). \qquad (3.7)$$

Substituting (3.7) into conditions 1 and 3 of assertion 3.1 and obtain:

1. Limit $\lim\limits_{n\to\infty} \dfrac{\int_2^n \frac{tf'(t)}{\log(t)}dt}{\frac{nf(n)}{\log(n)}} \neq 1$.

3. The limit $\lim\limits_{n\to\infty} \int_2^n \dfrac{tf'(t)}{\log(t)}dt = \pm\infty$.

Let's look at an example $\sum\limits_{p \leq n} \dfrac{1}{p}$. In this case: $f(n) = \dfrac{1}{n}$.

1. $\lim\limits_{n\to\infty} \dfrac{-\int_2^n \dfrac{dt}{t^3 \log(t)}}{\dfrac{n}{n\log(n)}} = -\lim\limits_{n\to\infty} \dfrac{n\log(n)}{n^4 \log(n)} = 0 \neq 1$.

2. $f(n) = \dfrac{1}{n}$ is a monotonic function and $f'(n) = -\dfrac{1}{n^2} \neq 0$.

3. $-\lim\limits_{n\to\infty} \int_2^n \dfrac{tdt}{t^2 \log(t)} = -\lim\limits_{n\to\infty}(\log\log(n)) = -\infty$.

Thus, in this case, all conditions of assertion 3.1 and corollary 3.2 are satisfied; therefore, sufficient conditions for the existence of asymptotic (2.5) are satisfied.

Similarly to the asymptotic formula (2.5), one can verify the fulfillment of the conditions of statement 3.1 and corollary 3.2 for the asymptotic formulas of other examples in chapter 2.

The conditions of assertion 3.1 and corollary 3.2 are sufficient for the indicated asymptotic to hold.

Let's consider a special case of assertion 3.1 and corollary 3.2, when the function is monotonically increasing and tends to infinity.

Assertion 3.3

Let $\lim\limits_{n\to\infty} f(n) = \infty$, $f'(n)$ is a continuous function, $f'(n) > 0$ and $\lim\limits_{n\to\infty} \dfrac{f(n)}{nf'(n)} \neq 0$. Then all conditions of assertion 3.1 and corollary 3.2 are satisfied.

Proof

Let's check the 1st condition:

$$\lim_{n\to\infty} \dfrac{\int_2^n \dfrac{tf'(t)}{\log(t)}dt}{\dfrac{nf(n)}{\log(n)}} = \lim_{n\to\infty} \dfrac{\dfrac{nf'(n)}{\log(n)}}{(\dfrac{nf(n)}{\log(n)})'} = \lim_{n\to\infty} \dfrac{\dfrac{nf'(n)}{\log(n)}}{\dfrac{nf'(n)}{\log(n)} + (\dfrac{n}{\log(n)})'f(n)} = \lim_{n\to\infty} \dfrac{1}{1 + \dfrac{f(n)}{nf'(n)}} \neq 1$$

Let's check the 2nd condition:

$f(n)$ is a monotone function and $f(n) \neq 0$ by the condition of assertion 3.3.

Let's check the 3rd condition:

$$\int_2^n \frac{tf'(t)}{\log(t)} dt \geq \int_2^n \frac{2f'(t)}{\log(2)} dt = \frac{2}{\log(2)} \int_2^n f'(t) dt = \frac{2}{\log(2)} (f(n) - f(2)) \quad - \quad \text{increases}$$

indefinitely, like $f(n)$.

Assertion 3.3 holds, for example, for the function $\sum_{p \leq n} p^m$, where $m > -1$.

Assertion 3.4

Let $f'$ is a continuous function and the above sufficient conditions are satisfied, then the following asymptotic equality holds:

$$\int_2^n \frac{f(t)dt}{\log(t)} = \left( \frac{f(n)n}{\log(n)} - \frac{2f(2)}{\log(2)} - \int_2^n \frac{tf'(t)dt}{\log(t)} \right)(1 + o(1)). \qquad (3.8)$$

Proof

Having in mind:

$$J_1(t) = \int_2^t \frac{du}{\log(u)} = \frac{t}{\log(t)} (1 + o(1)).$$

Integrating by parts, we get:

$$\int_2^n \frac{f(t)dt}{\log(t)} = f(t)J_1(t) \Big|_2^n - \int_2^n J_1(t) f'(t) dt = \left( \frac{f(n)n}{\log(n)} - \frac{2f(2)}{\log(2)} - \int_2^n \frac{tf'(t)dt}{\log(t)} \right)(1 + o(1)),$$

which corresponds to (3.8).

Let's consider an example of using assertion 3.4.

Let us prove the asymptotic equality:

$$\sum_{p \leq n} p^m = \frac{n^{m+1}}{(m+1)\log(n)} (1 + o(1)), \qquad (3.9)$$

for values $m > -1$.

Based on assertion 3.4:

$$J = \int_2^n \frac{t^m dt}{\log(t)} = \frac{n^{m+1}}{\log(n)} - \frac{2^{m+1}}{\log(2)} - mJ. \qquad (3.10)$$

Having in mind (3.10) we get:

$$J = \int_2^n \frac{t^m dt}{\log(t)} = \frac{1}{m+1}\left(\frac{n^{m+1}}{\log(n)} - \frac{2^{m+1}}{\log(2)}\right). \qquad (3.11)$$

Based on (3.11), we obtain the asymptotic:

$$\sum_{p \leq n} p^m = \sum_{k=2}^n \frac{k^m}{\log(k)}(1+o(1)) = J(1+o(1)) = \frac{n^{m+1}}{(m+1)\log(n)}(1+o(1)),$$

which corresponds to (3.9).

Asymptotic equality (3.9) corresponds to [3].

Let us now consider a necessary condition for the fulfillment of the indicated asymptotic.

Below we will use the notation of assertion 3.1.

Assertion 3.5

Let:

$$\lim_{n \to \infty} \frac{\sum_{k=1}^n a_k f(k)}{\sum_{k=1}^n b_k f(k)} = 1$$

Then at $p \to \infty$ ($p$ is a prime number):

$$\left|\frac{f(p)}{\sum_{k=1}^p b_k f(k)}\right| \to 0. \qquad (3.12)$$

Proof

If:

$$\lim_{n\to\infty} \frac{\sum_{k=1}^{n} a_k f(k)}{\sum_{k=1}^{n} b_k f(k)} = 1.$$

Then it runs:

$$\lim_{n\to\infty} \left| \frac{\sum_{k=1}^{n} a_k f(k)}{\sum_{k=1}^{n} b_k f(k)} - \frac{\sum_{k=1}^{n-1} a_k f(k)}{\sum_{k=1}^{n-1} b_k f(k)} \right| = 0. \qquad (3.13)$$

Take $n = p$, where $p$ is a prime number. Then, having in mind (3.13):

$$\lim_{p\to\infty} \left| \frac{\sum_{k=1}^{p} a_k f(k)}{\sum_{k=1}^{p} b_k f(k)} - \frac{\sum_{k=1}^{p-1} a_k f(k)}{\sum_{k=1}^{p-1} b_k f(k)} \right| = 0 \qquad (3.14)$$

Let us denote for $k = p$ by value $a_k f(k) = a_p f(p)$.

We transform (3.14) and taking into account that $a_p = 1$, we get:

$$\left| \frac{\sum_{k=1}^{p} a_k f(k)}{\sum_{k=1}^{p} b_k f(k)} - \frac{\sum_{k=1}^{p-1} a_k f(k)}{\sum_{k=1}^{p-1} b_k f(k)} \right| = \left| \frac{a_p f(p) \sum_{k=1}^{p-1} b_k f(k) - b_p f(p) \sum_{k=1}^{p-1} a_k f(k)}{\sum_{k=1}^{p} b_k f(k) \sum_{k=1}^{p-1} b_k f(k)} \right| = |f(p)| \left| \frac{1 - b_p \frac{\sum_{k=1}^{p-1} a_k f(k)}{\sum_{k=1}^{p-1} b_k f(k)}}{\sum_{k=1}^{p} b_k f(k)} \right|$$

Since for $p \to \infty$ the value $b_p \to 0$ and $\lim_{n\to\infty} \frac{\sum_{k=1}^{n} a_k f(k)}{\sum_{k=1}^{n} b_k f(k)} = 1$, then it is executed:

$$\left| \frac{f(p)}{\sum_{k=1}^{p} b_k f(k)} \right| \to 0 ,$$

which corresponds to (3.12).

For example, the specified prerequisite is not met for the function $f(p) = 2^p$ because

$$\sum_{k=1}^{p} b_k f(k) = \sum_{k=1}^{p} \frac{2^k}{\log(k)} = \frac{2^{p+1}}{\log p}.$$

## 4. CONCLUSION AND SUGGESTIONS FOR FURTHER WORK

The next article will continue to study the behavior of some sums.

## 5. ACKNOWLEDGEMENTS

Thanks to everyone who has contributed to the discussion of this paper. I am grateful to everyone who expressed their suggestions and comments in the course of this work.